\documentclass[a4paper,12pt]{amsart}

\usepackage{amssymb}
\usepackage{curves, epic}

\theoremstyle{plain}
\newtheorem{thm}{Theorem}[section]
\newtheorem{prop}[thm]{Proposition}
\newtheorem{lemma}[thm]{Lemma}
\newtheorem{corol}[thm]{Corollary}
\newtheorem*{ml}{Main Lemma}

\theoremstyle{remark}
\newtheorem*{remark}{Remark}

\numberwithin{equation}{section}

\begin{document}

\title[Fixed points of group actions]{Fixed points of discrete nilpotent\\ group actions on S$^2$}

\author{S. Druck} 
\thanks{The first and the third authors were partially supported by FINEP and FAPERJ of Brazil}

\author{F. Fang}

\thanks{The second author was partially supported by CNPq of Brazil; RFDP, 
Qiu-Shi Foundation and NSFC 19925104 of China}

\author{S. Firmo}



\date{September 1, 2001}

\address{S. Druck and S. Firmo\\ Universidade Federal Fluminense\\Instituto de Matem\'atica\\
Rua M\'ario Santos Braga s/n\\
24020-140 Niter\'oi, RJ\\Brazil}
\address{Fuquan Fang \\Nankai Institute of Mathematics, Tianjin 300071, P. R. China
\\Current address: Universidade Federal Fluminense\\Instituto de Matem\'atica\\
Rua M\'ario Santos Braga s/n\\
24020-140 Niter\'oi, RJ\\Brazil}

\email{druck@mat.uff.br}
\email{fang@mat.uff.br \\ffang@sun.nankai.edu.cn}
\email{firmo@mat.uff.br}

\maketitle

\begin{abstract}
We prove that for each integer \,$k\ge 2$\, there is an open neighborhood \,$\mathcal{V}_k$\, of
the identity map of the $2$-sphere $S^2$, in \,$C^1$ topology such that: if \,$G$\, is a 
nilpotent subgroup of \,$\textrm{Diff}^{\hskip1pt 1}(S^2)$\, with length \,$k$\, of nilpotency, generated by elements
in
\,$\mathcal{V}_k$\,, then the natural
\,$G$-action on 
\,$S^2$\, has nonempty fixed point set. Moreover,  the \,$G$-action  has at least two fixed points 
if the action has a finite nontrivial orbit.
\end{abstract}

\medskip
\medskip

\thispagestyle{empty}


\section{Introduction}
\vskip2mm


The classical Poincar\'e Theorem \cite{po} asserts that a \,$C^1$ vector field on a closed surface \,$\Sigma$\, 
with nonzero Euler characteristic has a singularity. Another way to phrase this conclusion is to
say that the flow tangent to the vector field must have a stationary point. 
In  \cite{l1,l2,l3} Lima  proved that pairwisely commuting vector fields on the surface \,$\Sigma$\, 
have a common singularity. This result has been generalized by Plante \cite{p1} for connected nilpotent 
Lie group actions  on \,$\Sigma$\,. The same result does not hold for solvable Lie group actions, as 
pointed out in the work of Lima.

From the Lefschetz Fixed Point Theorem any  diffeomorphism of the surface \,$\Sigma$\, homotopic to
the identity map has a fixed point. This strenghthened the Poincar\'e Theorem. A discrete analogue of 
Lima's Theorem was proved by Bonatti \cite{b1, b2}  asserting that commuting diffeomorphisms of  
\,$\Sigma$\, have a common fixed point, provided they are sufficiently $C^1$-close to the identity map.

In this paper we will prove a fixed point theorem for discrete nilpotent group actions on the $2$-sphere,
extending Bonatti's Theorem \cite{b1}. This may be considered as a discrete version of Plante's Theorem for
$S^2$.

\medskip

For a group \,$G$\,, its lower central series is         
$${G}={G}_{(0)}\supset{G}_{(1)}\supset\cdots
\supset{G}_{(k)}\supset{G}_{(k+1)}\supset\cdots$$
where \,${G}_{(k+1)}=[{G}\,,{G}_{(k)}]$\,. 
The first \,$k\geq 1$\, such that \,$G_{(k)}$\, is trivial is called the \,\emph{nilpotency length}\, of \,$G$.
If such \,$k$\, exists then \,$G$\, is called a \,\emph{nilpotent group}\, with nilpotency length \,$k$, or shortly,
a \,\emph{$k$-nilpotent group}.

The group of \,$C^1$\, diffeomorphisms of \,$S^2$\, endowed with the \,$C^1$ topology is denoted by 
\,$\textrm{\rm Diff}^{\hskip1pt 1}(S^2)$\,.  
For each subset \,$\mathcal{H}$\, of \,$\textrm{\rm Diff}^{\hskip1pt 1}(S^2)$\,, we denote by \,$Fix(\mathcal{H})$\,
the fixed point set of \,$\mathcal{H}$, i.e. 
$$Fix(\mathcal{H})=\{x\in S^2 \ ; \ f(x)=x \ , \ \forall f\in \mathcal{H} \}.$$ 
When \,$\mathcal{H}$\, is a subgroup of \,$\textrm{\rm Diff}^{\hskip1pt 1}(S^2)$\,, the fixed points of 
\,$\mathcal{H}$\, will be
called the fixed points of the \emph{natural
\,$\mathcal{H}$-action on
\,$S^2$} (shortly \,\emph{$\mathcal{H}$-action}).

\medskip

Our main results are as follows:

\medskip

\begin{thm}\label{teorcentral}
There exists a decreasing nested sequence \,$\{\mathcal{V}_k\}_{k\geq1}$  
of open neighborhoods of the identity map,  in \,$\textrm{\rm Diff}^{\hskip1pt 1}(S^2)$, such that\,$:$ if
\,$G$\, is a \,$k$-nilpotent subgroup of \,$\textrm{\rm Diff}^{\hskip1pt 1}(S^2)$\, finitely  generated by  elements
in 
\,$\mathcal{V}_{k}$\, then the \,$G$-action   has  a fixed point.
\end{thm}

\medskip

We refer to Section 5 for the precise definition of the open neighborhood \,$\mathcal{V}_k$\,, which is independent
of the number of generators for the group \,$G$.

Theorem \ref{teorcentral}  also holds for a \,$k$-nilpotent group \,$G$\, generated
by any family \,$\Gamma$\, of elements in \,$\mathcal{V}_{k}$\, for the following reasoning (c.f. Lima \cite{l3}).
Since  \,$G$\, is a \,$k$-nilpotent group then  any subgroup of \,$G$\, generated by a finite subset   of
\,$\Gamma$\,  is also \,nilpotent  with nilpotency length at most \,$k$.
By Theorem 1.1 the fixed point set  \,$Fix(\gamma)$\, is nonempty for any finite subset \,$\gamma$\, of \,$\Gamma$.
Therefore the family of closed subsets 
\,$\big(Fix(f)\big)_{f\in \Gamma}$\, of \,$S^2$\, has the ``finite intersection property'' which
implies that \,$Fix(\Gamma)\neq \emptyset$\,. Consequently we get

\medskip

\begin{corol}\label{corolinf}
Let \,$\mathcal{V}_k$\, be given by  Theorem \ref{teorcentral}. If \,$G\subset\textrm{\rm Diff}^{\hskip1pt
1}(S^2)$\, is a
\,$k$-nilpotent subgroup  generated by elements in \,$\mathcal{V}_k$\, then  the \,$G$-action   has a
 fixed point.
\end{corol}

\medskip

Theorem \ref{teorcentral} is just Bonatti's Theorem when \,$k=1$.

If the action in Theorem \ref{teorcentral} has a finite but \,{\it nontrivial}\, orbit, i.e. the orbit contains
at least two points, we have the following strengthened theorem. 

\medskip

\begin{thm}\label{teore3}
Let \,$G$\, and \,$\mathcal{V}_k$\, be as in Theorem \ref{teorcentral}. If the \,$G$-action has a finite  
nontrivial orbit then it has at least two  fixed points.
\end{thm}

\medskip

As an immediate corollary we have

\medskip

\begin{corol}\label{corol2fix}
Let \,$G$\,  be as in Theorem \ref{teorcentral} and let \,$Z(G)$ be the center of \,$G$. If there exists an element 
\,$h \in Z(G)$\, such that the number of fixed points of \,$h$\, satisfies 
\,$2\le \#\big({Fix}(h)\big)<\infty$\,  then the  
\,$G$-action  has at least two  fixed points.
\end{corol}

The above theorem should be compared with a result of Ghys \cite{gh1} which asserts that an analytic action of a
nilpotent group on \,$S^2$\, has always a finite orbit.

\medskip

Using the universal covering \,$S^2  \rightarrow \mathbb{RP}^2$\, it is easy to see that Theorems \ref{teorcentral}
and \ref{teore3}, and Corollaries \ref{corolinf} and \ref{corol2fix}  hold identically replacing \,$S^2$\, by
\,$\mathbb{RP}^2$.

\medskip

A typical application of Theorem \ref{teorcentral} to foliation theory is the following.

\medskip

\begin{thm}  
Let \,$F$\, be a closed manifold with nilpotent fundamental group. Then every foliation
\,$C^1$-close to the trivial foliation defined by the fibration \,$F\rightarrow F\times \Sigma^2 \rightarrow
\Sigma^2$\, must have a compact leaf close to some fiber, where \,$\Sigma^2=S^2$ or \,$\mathbb{RP}^2$.
\end{thm}

\medskip

Bonatti's Theorem was generalized by Handel \cite{ha} to commuting homeomorphisms of \,$S^2$\, which 
are \,$C^0$-close to the identity map. In view of  Handel and  Plante's Theorems, it is natural to
ask the following questions:

\begin{itemize}
\item [(i)] 
Do our theorems hold for nilpotent subgroups of homeomorphisms of \,$S^2$?

\item[(ii)] 
Does Plante's Theorem hold for discrete nilpotent group actions on surfaces with higher genus 
(at least $2$)?
\end{itemize}

In a forthcoming paper \cite{dff}  we will deal with the latter question, corresponding to Bonatti's
work [2].  
\medskip

We conclude the introduction with the following remark.

\medskip

An elementary result of Plante and Thurston \cite{E} asserts that finite dimensional connected nilpotent Lie 
subgroup of \,$\textrm{\rm Diff}^{\hskip1pt 1}(S^2)$\, is metabelian, i.e. its 
commutator subgroup is abelian. By Ghys \cite{gh1} any nilpotent subgroup of  
the real analytic diffeomorphisms group of \,$S^2$\,
is also metabelian. It is not known if the same result holds for $C^1$ diffeomorphisms.
On the other hand, there are examples (c.f \cite{E}) of connected nilpotent metabelian 
Lie subgroups of \,$\textrm{\rm Diff}^{\hskip1pt 1}(\mathbb{R}^2)$\, with arbitrarily large 
nilpotency length. 

\medskip
{\it Acknowledgement}: The paper was finished during the second author's visit to the Universidade
Federal Fluminense in Brasil. He would like to thank the institution for the hospitality.


\vskip40pt

\section{An algebraic proposition}
\vskip2mm


In this section we prove an algebraic proposition concerning with the choice of a set of generators for the
commutator subgroup of a nilpotent  group. First we need the following two lemmas.

\medskip

\begin{lemma}\label{homonilpot}
Let     \,$G$\, be a  \,$(k+1)$-nilpotent group where \,$k\geq1$. Then        
\begin{itemize}
\item [(1)]
$[f,h_1h_2]=[f,h_1] \, [f,h_2]$\, for all   \,$f\in G$\, and  \,$h_1,h_2\in G_{(k-1)}\,;$

\item[(2)]
$[f_1f_2,h]=[f_1,h] \, [f_2,h]$\, for all  \,$f_1,f_2\in G$\, and $h\in G_{(k-1)}$\,.
\end{itemize}
In particular, \,$[f,h^{-1}]=[f,h]^{-1}=[f^{-1},h]$\,  whenever                   
 \,$f\in G$\, and    \,$h\in G_{(k-1)}$\,.
\end{lemma}

\begin{proof}
Note that \,$G_{(k)}$\, is a subgroup of the center of \,$G$\, since \,$G$\, is \,$(k+1)$-nilpotent. 

Let \,$f,f_1,f_2\in G$\, and   \,$h,h_1,h_2\in G_{(k-1)}$\,. By definition
$[f, h_1]$ and $[f, h_2]$ are elements of $G_{(k)}$ and so they are in the center of
$G$. Therefore,
\begin{align}
 \nonumber [f,h_1h_2] & =  f(h_1h_2)f^{-1}(h_1h_2)^{-1}=  fh_1 f^{-1} h^{-1}_1 h_1 f h_2 f^{-1} h^{-1}_2
h^{-1}_1        \\ \nonumber
      & =   [f,h_1] \, h_1 \, [f,h_2] \, h^{-1}_1 = [f,h_1] \, [f,h_2]\,.\nonumber
\end{align}
Similarly, 
\begin{align}
\nonumber
[f_1f_2,h] & =  (f_1f_2)h(f_1f_2)^{-1}h^{-1} = f_1 [f_2,h] h f^{-1}_1 h^{-1}  \\ \nonumber
           & =    f_1 h f^{-1}_1 h^{-1} [f_2,h] =  [f_1,h] \,  [f_2,h]\,.    \nonumber
\end{align}
The proof is finished.         
\end{proof}

\medskip

For each nonempty subset   \,$\mathcal{S}$\, of a group \,$G$\, 
we set
$$\mathcal{S}_{(0)}=\mathcal{S} \quad \mbox{and} \quad \mathcal{S}_{(i+1)}=\{[a,b] \ ;  
\ a\in \mathcal{S} \ \ \mbox{and} \ \
b\in\mathcal{S}_{(i)}\}\,, \ i\ge 0\,.$$

We shall denote by  \,$\langle \mathcal{H} \rangle$\, the  subgroup  generated by the subset \,$\mathcal{H}$\, of a
group.

\medskip

\begin{lemma}\label{geragrk} 
If \,$G$\, is a \,$k$-nilpotent group generated by \,$\mathcal{S}$\,                  
 then  we have   \,\mbox{$G_{(k-1)}=\langle \mathcal{S}_{(k-1)}\rangle$}.
\end{lemma}

\begin{proof}
We prove this lemma by induction on  \,$k$.

The   lemma is trivial for \,$k=1$\, since in this case we have 
$$G_{(0)}=G=\langle S \rangle = \langle S_{(0)} \rangle .$$

Suppose the lemma is true for some \, $k\geq 1$\,.  We need to prove it for \,$k+1$.

For this, let \,$G$\, be a  \,$(k+1)$-nilpotent group and let \,$\mathcal{S}$\, be a set of generators for \,$G$.
Consider the lower central series   
\begin{equation}\label{seqinde}
 G \varsupsetneq G_{(1)}\varsupsetneq \cdots  \varsupsetneq G_{(k-1)}  \varsupsetneq G_{(k)} \varsupsetneq
G_{(k+1)}=[G,G_{(k)}]=\{e\}.
\end{equation}

Since \,$G_{(k)}$\, is a subgroup of the center of \,$G$, we get  
$$\frac{G_{(i)}}{G_{(k)}} = \Big(\frac{G}{G_{(k)}}\Big)_{(i)} \quad \mbox{for all  } \quad 0\leq i \leq k\,.$$
By \ref{seqinde}  we have
\begin{equation} 
\frac{G}{G_{(k)}}  \varsupsetneq \Big(\frac{G}{G_{(k)}}\Big)_{(1)} \varsupsetneq \cdots  \varsupsetneq 
\Big(\frac{G}{G_{(k)}}\Big)_{(k-1)}
\varsupsetneq \Big(\frac{G}{G_{(k)}}\Big)_{(k)}
=\{e\}.\nonumber
\end{equation}
Therefore \,$G/G_{(k)}$\, is  a \,$k$-nilpotent group generated by the projection  \,$\
\overline{\!\mathcal{S}}$\,  of  \,$\mathcal{S}$\, onto \,$G_{}/G_{(k)}$\, by the quotient map.
By the induction assumption we have that  
$$\frac{G_{(k-1)}}{G_{(k)}} = \Big(\frac{G}{G_{(k)}}\Big)_{(k-1)} = \langle \ \overline{\!\mathcal{S}}_{(k-1)}
\rangle $$ where \,$\ \overline{\!\mathcal{S}}_{(k-1)}$\, denotes the projection of \,$\mathcal{S}_{(k-1)}$\, onto 
\,$G/G_{(k)}$\, by the quotient map. Therefore for each \,$h\in G_{(k-1)}$\, there exist
\begin{equation}
h_1,\ldots,h_n \in \mathcal{S}_{(k-1)} \quad \mbox{and} \quad \epsilon_1,\ldots,\epsilon_n\in\{1,-1\}\nonumber
\end{equation}
such that          
\begin{equation}
\overline{h} = 
 \big(\,\overline{h}_1\big)^{\epsilon_1} \ldots \big(\,\overline{h}_n\big)^{\epsilon_n}.\nonumber
\end{equation}
Consequently, there exists a \,$\xi \in G_{(k)}$\, such that       
\begin{equation}\label{prodgk}
{h} = 
{h}_1^{\epsilon_1} \ldots{h}_n^{\epsilon_n} \, \xi.
\end{equation}

Note that each \,$f\in G$\, can be written as 
\,$f=f_1^{\delta_1} \ldots f_m^{\delta_m}$\, where 
\,$f_1,\ldots,f_m \in \mathcal{S}$\, and  \,$\delta_1,\ldots, \delta_m\in\{1,-1\}$\,.
 From      \,$[G,G_{(k)}]=\{e\}$\,, 
the identity \ref{prodgk} and the Lemma \ref{homonilpot} we get that   
$$[f,h]=\prod_{\substack{1\leq i\leq m \\ 1\leq j\leq n}} [f_i,h_j]^{\delta_i \epsilon_j} [f_i,\xi]^{\delta_i}=
\prod_{\substack{1\leq i\leq m \\ 1\leq j\leq n}} [f_i,h_j]^{\delta_i \epsilon_j}\,. $$
Since     \,$[f_i,h_j]\in \mathcal{S}_{(k)}$\, 
this proves that    \,$G_{(k)}$\, is generated by \,$\mathcal{S}_{(k)}$\, and the proof is finished.
\end{proof}

\medskip

\begin{prop}\label{geragrgk}
If  \,$G$\, is a \,$k$-nilpotent group generated by   \,$\mathcal{S}$\,     
then    \,$G_{(1)}=\langle \mathcal{S}_{(1)},\ldots, \mathcal{S}_{(k)}\rangle .$
\end{prop}

\begin{proof}
We use induction on \,$k$\, once again in the proof.  

The case  \,$k=1$\, is trivial. The case \,$k=2$\, follows from Lemma \ref{geragrk}.

Suppose the lemma is true for some  \,$k\geq 2$\,.

Let \,$G$\, be a \,$(k+1)$-nilpotent group and let \,$\mathcal{S}$\, be
a set of generators for \,$G$. By the proof of the last lemma  we have  
$$\frac{G_{(i)}}{G_{(k)}} = \Big(\frac{G}{G_{(k)}}\Big)_{(i)} \quad \mbox{for all} \quad 0\leq i \leq k$$
and \,$G/G_{(k)}$\, is  \,$k$-nilpotent.

By the induction assumption we get  
\begin{equation}
\frac{G_{(1)}}{G_{(k)}} = \langle \,\, \overline{\!\mathcal{S}}_{(1)},\ldots,\,\overline{\!\mathcal{S}}_{(k-1)}
\rangle .\nonumber
\end{equation}
Consequently, for each \,$h\in G_{(1)}$\, its image  \,$\overline{h}$\, in the quotient group
\,${G_{(1)}}/{G_{(k)}}$\, can be written as  
\begin{equation}
\overline{h}= \big(\,\overline{h}_1\big)^{\epsilon_1} \ldots \big(\,\overline{h}_n\big)^{\epsilon_n}
\nonumber
\end{equation}
where 
$$h_1,\ldots,h_n\in \bigcup_{i=1}^{k-1}  \mathcal{S}_{(i)}  \quad \mbox{and} \quad 
\epsilon_1,\ldots,\epsilon_n\in\{1,-1\}.$$
It is easy to see that 
there exists a \,$\xi\in G_{(k)}$\, so that       
\begin{equation}\label{decompggk}
h=  {h}_1^{\epsilon_1} \ldots h_n^{\epsilon_n} \, \xi.\nonumber
\end{equation}
On the other hand,   by Lemma \ref{geragrk} there exist 
$$\xi_{1},\ldots,\xi_{m}\in \mathcal{S}_{(k)} \quad \mbox{and} \quad \delta_{1},\ldots,\delta_{m}\in\{1,-1\}$$ 
so that         
\begin{equation}\label{decompgk} 
\xi = {\xi}_{1}^{\delta_{1}} \ldots \xi_{m}^{\delta_{m}}\,.\nonumber
\end{equation}
Therefore \,$G_{(1)}$\, is generated by \,$\bigcup_{1\leq i\leq k}\mathcal{S}_{(i)}$\, 
and the proof is finished.
\end{proof}


\vskip40pt

\section{Invariance and recurrence}
\vskip2mm


If  \,$f,g$\, are commuting diffeomorphisms then  \,${Fix}(g)$\,
is \,$f$-invariant. In this section we show that the fixed point set of the commutator 
subgroup, $Fix(G_{(1)})$\, with \,$G\subset \textrm{Diff}^{\hskip1pt 1}(S^2)$,  has invariance and recurrence
properties even without the assumption of commutativity. These  properties play an important role in this paper.

\medskip

Let \,$f\in\textrm{Diff}^{\hskip1pt 1}(S^2)$\,. The \,\emph{positive semi-orbit}\, of a point \,$p\in S^2$\, is the
set \,$\mathcal{O}_p^+(f)=\{f^i(p) \ ; \ i\in\mathbb{N}\}$. Its closure will be denoted by
\,$\overline{\mathcal{O}_p^+(f)}$\,.

We say that  \,$p\in S^2$\, is  \,{\it  $\omega$-recurrent}\, for   
  \,$f$\, if  
\,$p$\, is the limit of some subsequence of \,$\big(f^{n}(p)\big)_{n\in\mathbb{N}}$.

\medskip

\begin{prop}\label{invariorbit}
Let  \,$G\subset\textrm{\rm Diff}^{\hskip1pt 1}(S^2)$\, be a subgroup  and let 
$$f,f_1,\ldots,f_n\in G .$$ Then
\,$Fix(G_{(1)}\,,f_1,\ldots,f_n)$\, is  \,$f$-invariant. 
Moreover, it contains \,$\omega$-recurrent points for \,$f$\, lying in  
\,$\overline{\mathcal{O}_p^+(f)}$, for all \,$p\in Fix(G_{(1)}\,,f_1,\ldots,f_n)$. 
\end{prop}

\begin{proof}
The first assertion is an immediate consequence of next lemma.
To prove the second assertion let 
\,$p\in Fix(G_{(1)},f_1,\ldots,f_n)$\,. Since
\,$Fix(G_{(1)},f_1,\ldots,f_n)$\, is \,$f$-invariant it follows that  
$$\overline{\mathcal{O}^{+}_{p}(f)} \subset Fix(G_{(1)},f_1,\ldots,f_n)\,.$$
 By Zorn's Lemma  \,$f$\, has a minimal set in 
\,$\overline{\mathcal{O}^{+}_{p}(f)}$\,. The points of \,$\overline{\mathcal{O}^{+}_{p}(f)}$\, in a minimal set of
\,$f$\, are 
\,$\omega$-recurrent points for \,$f$. This completes the proof.
\end{proof}

\medskip

\begin{lemma}
Let \,$G$\, be as in Proposition \ref{invariorbit}. Then
$Fix(G_{(i)},g)$\, is \,$f$-invariant for all \,$f\in G$\, and \,$g\in G_{(i-1)}$\, where \,$i\geq1$.
\end{lemma}

\begin{proof}
Let  \,$p\in Fix(G_{(i)},g)$\,. As \,$p\in Fix(G_{(i)})$\, we have that                      
\,$p\in Fix(G_{(i+1)})$\,.  
For  \,$h\in G_{(i)} \,\, , \, g\in G_{(i-1)}$\, and \,$f\in G$\, we have
$$[f^{-1},h]\in G_{(i+1)} \ \ \mbox{and} \ \ [f^{-1},g]\in G_{(i)} \,.$$
Thus
\begin{eqnarray}
p=[f^{-1},h](p)=f^{-1}hfh^{-1}(p)=f^{-1}hf(p)\nonumber
\end{eqnarray}
and
\begin{eqnarray}
p=[f^{-1},g](p)=f^{-1}g f g^{-1} (p)=f^{-1} g f(p)\,.\nonumber
\end{eqnarray}
Therefore 
$$hf(p)=f(p) \ \ \mbox{and} \ \  gf(p)=f(p)$$ 
which prove that \,$Fix(G_{(i)},g)$\, is \,$f$-invariant.
\end{proof}

\medskip

\begin{remark}
The results of this section hold for any group of homemorphisms of an \,$n$-dimensional manifold \,$M$\, where
\,$M$\, needs to be compact for the second assertion in Proposition \ref{invariorbit}.
\end{remark}


\vskip40pt

\section{A character curve}
\vskip2mm


In this section we generalize Lemma 4.1 of \cite{b1} by proving it without the commutativity hypothesis on \,$f$\,
and \,$g$\,  (c.f. Lemmas \ref{curvagama3} and \ref{curvagama1} below). This generalization will be an important
tool in the proofs of our theorems.

First, let us fix some notations and definitions. 

Let \,$S^2$\, denote the unit \,$2$-sphere in \,$\mathbb{R}^3$\, endowed with the standard norm,
denoted by   \,$||\ ||$\,.
If \,$\varphi:S^2\rightarrow\mathbb{R}^3$\, is a \,$C^1$ map, we define 
$$||\varphi||_1=\sup_{x\in S^2}\Big\{||\varphi(x)||+\sup_{v\in T_{x}S^2\,;\,||v||=1}||D\varphi(x).v||\Big\}\,.$$

For \,$a,b\in S^2$\, where \,$a\neq -b$\,, let  \,$[a,b]$\, denote the 
oriented minimal geodesic segment from \,$a$\, to \,$b$\,, and let \,$d(a,b)$\, denote its length.   

For  \,$f\in\mbox{Diff}^{\hskip1pt 1}(S^2)$\, without antipodal points and \,$p\in S^2$\, let
$$\gamma^{\,p}_{f}:[0,\mu)\rightarrow S^2$$ 
be the curve obtained  joining the oriented minimal geodesic segments
\,$[f^i(p)\,,f^{i+1}(p)]$\, where \,$i\in\mathbb{N}$\, and \,$0<\mu\leq\infty$\,. 
Note that \,$\mu =\infty$\, if
\,$p$\,  is a nonfixed \,$\omega$-recurrent point for \,$f$.

Now let us fix  the following neighborhood of the identity map of \,$S^2$:  
$$\mathcal{V}_1=\big\{f\in\textrm{Diff}^{\hskip1pt 1}(S^2) \ ; \ ||f-\textrm{Id}||_1<\mbox{$\frac{1}{60}$}\big\}.$$

We recall the following two results of Bonatti  
\cite{b1} concerning  the curve \,$\gamma^{\,p}_{f}$\, which will be used in the proofs of our results.

\medskip

\begin{prop}\label{propbona1}{\rm\textbf{(Bonatti)}}
Let \,$f\in\mathcal{V}_1$.
\begin{itemize}
\item [(1)]
If \,$p\in S^2-Fix(f)$\, then \,$f$\, does not have any fixed point in the open
ball  \,$\mbox{\rm B}(p\,,4d(p\,,f(p)))$.\\
In particular, \,$f$\, does not have fixed points along \,$\gamma^{\,p}_{f}$.

\item[(2)]
If \,$p\in S^2-Fix(f)$\, is an \,$\omega$-recurrent point for \,$f$\, then there
exists a simple closed curve \,$\gamma$\, contained in \,$\gamma_{f}^{\,p}$.
\end{itemize}
\end{prop}

\medskip

We will call the simple closed curve \,$\gamma$\, obtained in the above proposition the 
\,\emph{character curve}\, of \,$f$\, at \,$p$\, where \,$p$\, is an \,$\omega$-recurrent point for \,$f$.

\noindent
\hfil
\setlength{\unitlength}{1mm}
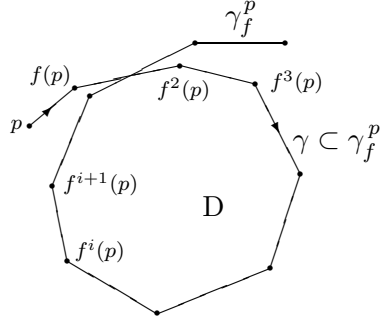
\begin{picture}(48,43)(-7,2)
\put(0,0){\drawline(-5,28)(1,33)(15,36)(25,33.5)(31,21.5)(27,9)(12,3)(0,10)(-2,20)(3,32)(17,39)(29,39)}

\put(-2,30.57){\vector(1,1){0}}
\put(-5,28){\circle*{.8}}
\put(-7.5,27.6){\tiny$p$}
\put(1,33){\circle*{.8}}
\put(-5,34){\tiny$f(p)$}
\put(12,32){\tiny$f^2(p)$}
\put(27,33){\tiny$f^3(p)$}
\put(15,36){\circle*{.8}}
\put(25,33.5){\circle*{.8}}
\put(31,21.5){\circle*{.8}}
\put(30,25){\small$\gamma\subset \gamma^{\,p}_{f}$}
\put(27,9){\circle*{.8}}
\put(1,10.5){\tiny$f^i(p)$}
\put(12,3){\circle*{.8}}
\put(-.5,19.4){\tiny$f^{i+1}(p)$}
\put(0,10){\circle*{.8}}
\put(-2,20){\circle*{.8}}
\put(18,16){\small D}
\put(3,32){\circle*{.8}}
\put(17,39){\circle*{.8}}
\put(29,39){\circle*{.8}}
\put(21,41.5){\small$\gamma^{\hskip1pt p}_{f}$}
\put(28.23,27){\vector(1,-2){0}}
\end{picture}

\medskip

\begin{prop}\label{propbona2}{\rm\textbf{(Bonatti)}}
Let   \,$f,h_1\,,\ldots\,,h_n\in\mathcal{V}_1$\, be commuting diffeomorphisms and \,$n\in\mathbb{Z}^+$\,.
Let \,$p\in  Fix(h_1\,,\ldots\,,h_{n})- Fix(f)$\, be an \,$\omega$-recurrent point for 
\,$f$\, and \,$\gamma\subset \gamma^{\,p}_{f}$\, its character
curve at \,$p$.  If   
\,$\mbox{\rm D}$\, is a disk enclosed by   \,$\gamma$\, then
\,$f,h_1\,,\ldots,h_n$\, have a common fixed point in the interior of 
\,$\mbox{\rm D}$. 
\end{prop}

\medskip

Let \,$G\subset \textrm{Diff}^{\hskip1pt 1}(S^2)$\, be a  subgroup and 
$$f,g,h_1,\ldots,h_n\in  \mathcal{V}_1 \cap G\,.$$

Suppose there exist two points \,$p\,,q\in S^2$ satisfying
\begin{itemize}
\item []
$p\in {Fix}(G_{(1)},h_1\,,\ldots,h_n\,,g)- {Fix}(f)$\,;
\item []
$q \in {Fix}(G_{(1)},h_1\,,\ldots,h_n\,,f)- {Fix}(g)$\,.
\end{itemize}

\medskip

With these assumptions we have the following three results.

\medskip

\begin{lemma}\label{curvagama3} 
The curves  \,$\gamma^{\,p}_{f}$\, and \,$\gamma^{\,q}_{\,g}$\, 
are disjoint.
\end{lemma}

\begin{proof}
Suppose not. Then there exist \,$i,j\in\mathbb{N}$\, so that 
$$[g^i(q)\,,g^{i+1}(q)]\cap[f^j(p)\,,f^{j+1}(p)]\neq\emptyset\,.$$
By the triangle inequality we have  
\begin{align}
 d\big(g^i(q),f^j(p)\big) & \leq   
d\big(g^i(q),g^{i+1}(q)\big)+d\big(f^j(p),f^{j+1}(p)\big) 
\nonumber\\
& \leq 2\max\Big(d\big(g^i(q),g^{i+1}(q)\big) \,,\, d\big(f^j(p),f^{j+1}(p)\big)\Big)\,. 
\nonumber
\end{align}

Consequently,
\begin{itemize}
\item [--]
either \,$f^j(p)$\, is  in the ball 
 \,$\textrm{B}\big(g^i(q),3d(g^i(q),g^{i+1}(q))\big)$\,, which is impossible
by Proposition \ref{invariorbit} and  Proposition \ref{propbona1}, since        
\,$f^j(p)$\, is a fixed point of \,$g$\,;

\item[--]
or \,$g^i(q)$\, is in the ball \,$\textrm{B}(f^j(p),3d(f^j(p),f^{j+1}(p)))$\,, 
which is impossible by  Proposition \ref{invariorbit} and Proposition \ref{propbona1} once again,
since \,$g^i(q)$\, is a fixed point of  \,$f$\,.
\end{itemize}
This shows that                 
\,$\gamma^{\,p}_{f}$\, and \,$\gamma^{\,q}_{\,g}$\, are disjoint curves.   
\end{proof}

\medskip

\begin{corol}\label{curvagama2}
Suppose that \,$q$\, is an \,$\omega$-recurrent point for \,$g$\, and \,$p\in \textrm{\rm Int(D)}$\,, where
\,$\textrm{\rm D}$\, is a disk enclosed by the character curve \,$\gamma\subset\gamma^{\,q}_{\,g}$. 
Then \,$\gamma^{\,p}_{f}\subset \mbox{\rm Int(D)}$\,. 
In particular, \,$\overline{\mathcal{O}^+_{p}(f)} \subset \mbox{\rm Int(D)}$\,.

\end{corol}

\begin{proof}
By Lemma \ref{curvagama3} we have that \,$\gamma^{\,p}_{f}\subset \mbox{\rm Int(D)}$\, since
\,$p\in\textrm{Int(D)}$\, and the curves \,$\gamma^{\,p}_{f}$\, and \,$\gamma^{\,q}_{\,g}$\, 
are disjoint. On the other hand, since
$$\overline{\mathcal{O}^+_p(f)} \subset Fix(g)$$
and \,$g$\, does not have fixed points along \,$\gamma^{\,q}_{\,g}$\, then we get
$$\overline{\mathcal{O}^{+}_{p}(f)}\subset \textrm{\rm Int(D)}\,.$$ 
This
completes the proof.
\end{proof}

\medskip

\begin{lemma}\label{curvagama1}
Suppose  there exists an \,$r>0$\, such that 
$$d\big(f^i(p)\,,f^{i+1}(p)\big)\geq r \quad \mbox{and} 
\quad d\big(g^i(q)\,,g^{i+1}(q)\big)\geq r$$
for all \,$i\in\mathbb{N}$\,.
Then  \,$d\big(\gamma^{\,p}_{f}\,,\gamma^{\,q}_{\,g}\big)\geq r.$
\end{lemma}

\begin{proof}
Suppose not. Then there are two points  
$$a\in[g^i(q)\,,g^{i+1}(q)]  \ \ \mbox{and} 
\ \  b\in[f^j(p)\,,f^{j+1}(p))] \ \ \mbox{such that} \ \
d(a,b)< r.$$ 
Therefore
\begin{align}
d\big(g^i(q),f^j(p)\big) & \leq
d\big(g^i(q),g^{i+1}(q)\big)+r+d\big(f^j(p),f^{j+1}(p)\big) 
\nonumber\\
  & \leq  
3\max\Big(d\big(g^i(q),g^{i+1}(q)\big) \,,\, d\big(f^j(p),f^{j+1}(p)\big)\Big)\,. 
\nonumber
\end{align}

Now the same argument used in the proof of 
Lemma \ref{curvagama3}  applies to conclude the lemma. 
\end{proof}


\vskip40pt

\section{The main lemma}
\vskip2mm


For each integer  \,$k\geq 1$\, let \,$\mathcal{V}_k$\, be the following open neighborhood of the identity map of
\,$S^2$\, in the
\,$C^1$ topology:
$$\mathcal{V}_k=\bigg\{ f\in \textrm{\rm Diff}^{\hskip1pt 1}(S^2) \ ; \ 
||f-\textrm{Id}||_1\le \frac {1}{5^{\frac{(k-1)k}{2}}\cdot 60} \bigg\}.$$

\medskip

It is an elementary exercise to verify that

\medskip

\begin{prop}\label{propviz} 
If \,$f, g \in \mathcal{V}_1$\, then 
$$||[f,g]-\mbox{\rm Id}||_1\le 5 \max\big\{||f-\mbox{\rm Id}||_1 \, , ||g-\mbox{\rm Id}||_1\big\}.$$
Furthermore,
if \,$f_1, \ldots ,f_{k+1} \in \mathcal{V}_{k+1}$\, and \,$k\geq1$\, then 
$$[f_1, [f_2, \ldots , [f_{i}, f_{i+1}]\ldots]] \in  \mathcal{V}_{k} \ \ \mbox{for all} \ \ 1\leq i\leq k\,.$$
\end{prop}

\medskip

For the sake of simplicity we use \,$f_0$\, to denote the identity map.

\medskip

\begin{ml}\label{fixoinducao} 
Let \,$G\subset \textrm{\rm Diff}^{\hskip1pt 1}(S^2)$\, be a \,$k$-nilpotent subgroup finitely
generated by elements in \,$\mathcal{V}_k$\, and let
\,$f_1\,,\ldots,f_n\in \mathcal{V}_k\cap G$.  Let
$$p\in Fix(G_{(1)},f_0\,,\ldots\,,f_{n-1})-{Fix}(f_{n})$$ 
be an \,$\omega$-recurrent point for  
\,$f_{n}$\, and \,$\gamma\subset \gamma^{\,p}_{f_{n}}$\, be the character curve of $f_n$ at $p$.
If \,$\mbox{\rm D}$ is a disk enclosed by \,$\gamma$\, then
\,$G_{(1)}\,,f_1\,,\ldots\,,f_n$\, have a common fixed point in the interior of \,$\textrm{\rm D}$.
\end{ml}

\begin{proof}
The proof will be by induction on the nilpotency length of the group \,$G$.
When \,$k=1$\,, the group  \,$G$\, is abelian and the lemma is just
Proposition \ref{propbona2}.
Assume now that for some \,$k\geq1$\, the lemma is true whenever the nilpotency length of \,$G$\, is \,$l$,  for all 
\,$1\leq l\leq k$\, and for all \,$n\in\mathbb{Z}^+$. 
Suppose that \,$G$\, is a \,$(k+1)$-nilpotent group finitely generated by  elements in \,$\mathcal{V}_{k+1}$\,
as in the lemma. We now proceed by induction on the number \,$n$\, of diffeomorphisms 
\,$f_1\,,\ldots\,,f_n$.

Let 
\,$f_{1}\in\mathcal{V}_{k+1}\cap G$\, and let \,$p\in Fix(G_{(1)})-{Fix}(f_{1})$\,  
be an \,$\omega$-recurrent point for \,$f_{1}$\,. Let      \,$\gamma\subset
\gamma^{\,p}_{f_{1}}$\,  be   the character curve of \,$f_1$\, and let \,$\textrm{D}$\, be a disk enclosed by
\,$\gamma$. Joining the Propositions 
\ref{geragrgk} and \ref{propviz} we obtain a finite set \,$\{g_1,\ldots,g_m\}\subset\mathcal{V}_{k}$\, which
generates
\,$G_{(1)}$\,.
Let \,$H$\, be the
subgroup of \,$G$\, generated by \,$G_{(1)}$\, and \,$f_{1}$\,. Note that  
\,$H_{(1)}\subset G_{(2)}$\,. Therefore \,$H$\, is a nilpotent group with  nilpotency length  at
most \,$k$\, and 
$$p\in Fix(H_{(1)},g_1,\ldots,g_m)-Fix(f_{1})\,.$$
By the induction assumption there exists a point \,$q$\, in the interior of the disk 
\,$\textrm{D}$\, such that 
$$q\in Fix({H_{(1)},g_1,\ldots,g_m,f_{1}})=Fix (G_{(1)},f_1)\,.$$
This proves the case  \,$n=1$.

Assume now that the lemma holds for some \,$n\geq1$\,. In order  to prove it for \,$n+1$\, 
let \,$f_1,\ldots,f_{n+1}\in\mathcal{V}_{k+1}\cap G$, and let 
$$p\in Fix(G_{(1)},f_{1},\ldots,f_{n})-Fix(f_{n+1})$$
be an \,$\omega$-recurrent point for \,$f_{n+1}$\,. Let \,$\gamma\subset\gamma^{\,p}_{f_{n+1}}$\, be the
character curve of \,$\gamma^{\,p}_{f_{n+1}}$\, and let us fix a disk \,$\textrm{D}$\, enclosed by \,$\gamma$.

Suppose by contradiction that
\,$G_{(1)},f_{1},\ldots,f_{n},f_{n+1}$\,
 do not have common fixed points in the interior of \,$\textrm{D}$. The induction assumption on \,$n$\, asserts
that \,$G_{(1)},f_{0},\ldots,f_{n-1},f_{n+1}$\, have a  common fixed point \,${\tilde y}_0$\, in
\,$\textrm{Int(D)}$\,. In view of our contradiction assumption we have that 
$${\tilde y}_0 \in Fix(G_{(1)},f_{0},\ldots,f_{n-1},f_{n+1})-Fix(f_{n})\,.$$
Applying Corollary \ref{curvagama2} to the diffeomorphisms \,$f_{n}$\, and \,$f_{n+1}$\, we conclude that
\,$\overline{\mathcal{O}^+_{{\tilde y}_0}(f_{n})}\subset \textrm{Int(D)}\,.$
According to Proposition \ref{invariorbit} the map \,$f_n$\, has an \,$\omega$-recurrent point
\,$y_0\in\overline{\mathcal{O}^+_{{\tilde y}_0}(f_{n})} \subset Fix(G_{(1)},f_{0},\ldots,f_{n-1},f_{n+1})$. 
Observe that \,$y_0\in\textrm{Int(D)}$\, is a  nonfixed 
\,$\omega$-recurrent point for \,$f_n$\,, that is
$$y_0\in Fix(G_{(1)},f_{0},\ldots,f_{n-1},f_{n+1})-Fix(f_{n})$$
From  Proposition \ref{propbona1} and Corollary  \ref{curvagama2}  we obtain a simple closed
 curve \,$\gamma_0\subset \gamma_{f_{n}}^{y_0}$\, and a closed disk $\textrm {D}_0$ enclosed by 
$\gamma _0$ so that \,$\textrm{D}_0\subset \textrm{Int(D)}$\,.

\noindent
\hfil
\setlength{\unitlength}{1mm}
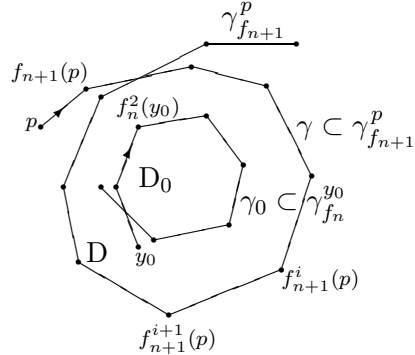
\begin{picture}(50,49)(-7,-4)
\put(0,0){\drawline(-5,28)(1,33)(15,36)(25,33.5)(31,21.5)(27,9)(12,3)(0,10)(-2,20)(3,32)(17,39)(29,39)}

\put(0,0){\drawline(8,12)(5,20)(8,28)(17,29.5)(21.9,23)(20,15)(10,13)(3,20)}
\put(5,29.8){\tiny$f^2_{n}(y_0)$}
\put(7.5,10){\tiny$y_0$}
\put(8,12){\circle*{.8}}
\put(5,20){\circle*{.8}}
\put(8,28){\circle*{.8}}
\put(17,29.5){\circle*{.8}}
\put(21.9,23){\circle*{.8}}
\put(20,15){\circle*{.8}}
\put(10,13){\circle*{.8}}
\put(3,20){\circle*{.8}}
\put(21.5,17){\small$\gamma_0\subset\gamma^{y_0}_{f_n}$}

\put(-2,30.57){\vector(1,1){0}}
\put(7,25.2){\vector(1,2){0}}

\put(-5,28){\circle*{.8}}
\put(-7,28){\tiny$p$}

\put(1,33){\circle*{.8}}
\put(-9,34.5){\tiny$f_{n+1}(p)$}

\put(15,36){\circle*{.8}}
\put(25,33.5){\circle*{.8}}

\put(31,21.5){\circle*{.8}}
\put(29,27){\small$\gamma\subset \gamma^{\hskip1pt p}_{f_{n+1}}$}

\put(27,9){\circle*{.8}}
\put(27,7){\tiny$f^i_{n+1}(p)$}

\put(12,3){\circle*{.8}}
\put(8,-1){\tiny$f^{i+1}_{n+1}(p)$}

\put(0,10){\circle*{.8}}
\put(-2,20){\circle*{.8}}
\put(1,10){\small D}
\put(8,20){\small D$_0$}

\put(3,32){\circle*{.8}}
\put(17,39){\circle*{.8}}
\put(29,39){\circle*{.8}}
\put(19,41.5){\small$\gamma^{\hskip1pt p}_{f_{n+1}}$}

\end{picture}

Once again  the induction assumption on \,$n$\, asserts that there exists a point
\,${\tilde x}_1\in \textrm{Int(D}_0)$\, such that
$${\tilde x}_1\in{Fix}(G_{(1)},f_{0}\,,\ldots,f_{n-1}\,,f_{n})-{Fix}(f_{n+1})\,.$$
Now we shall repeat the above argument.  Corollary \ref{curvagama2} applied to the diffeomorphisms \,$f_{n}$\, and
\,$f_{n+1}$\,  gives   
\,$\overline{\mathcal{O}^+_{{\tilde x}_1}(f_{n+1})}\subset  \textrm{Int(D}_0)$\, and Proposition \ref{invariorbit}
asserts that \,$f_{n+1}$\, has an  \,$\omega$-recurrent point  
$$x_1\in\overline{\mathcal{O}^+_{{\tilde
x}_1}(f_{n+1})} \subset {Fix}(G_{(1)},f_{0}\,,\ldots,f_{n-1}\,,f_{n}).$$
Thereby \,$x_1\in\textrm{Int(D}_0)$\,  is an \,$\omega$-recurrent point for \,$f_{n+1}$\, and
$$x_1 \in {Fix}(G_{(1)},f_{0}\,,\ldots,f_{n-1}\,,f_{n})-{Fix}(f_{n+1})\,.$$
Proposition \ref{propbona1} and Corollary  \ref{curvagama2}  give a simple closed curve
\,$\gamma_1\subset \gamma_{f_{n+1}}^{\,x_1}$\, and an enclosed disk  
\,$\textrm{D}_1\subset \textrm{Int(D}_0)$\,. 

Remember we are assuming that \,$G_{(1)},f_1,\ldots,f_n,f_{n+1}$\,  do not have common fixed points in the interior
of \,$\textrm{D}$\,. Consequently, they do not have common fixed points in  \,$\textrm{D}$\, since by Proposition 
\ref{propbona1}  the diffeomorphism \,$f_{n+1}$\, does not have fixed points along \,$\gamma$.

By the compactness of the sets
$$Fix(G_{(1)},f_0,\ldots,f_{n-1},f_{n})\cap \textrm{D} \quad \mbox{and} \quad 
Fix(G_{(1)},f_0,\ldots,f_{n-1},f_{n+1})\cap \textrm{D}$$ 
we get a constant \,$r>0$\, such that for all

\begin{itemize}
\item []
$x\in Fix(G_{(1)},f_{0},\ldots,f_{n-1},f_{n})\cap \textrm{D}$\, and

\item []
$y\in Fix(G_{(1)},f_{0},\ldots,f_{n-1},f_{n+1})\cap \textrm{D}$

\end{itemize}
the distance \,$d$\, satisfies 
\begin{eqnarray}\label{ineqinv}
d\big(x,f_{n+1}(x)\big)\geq r \quad \mbox{and}  \quad d\big(y,f_{n}(y)\big)\geq r. 
\end{eqnarray}
Since \,$\mathcal{O}^+_{y_0}(f_n)$\, and  \,$\mathcal{O}^+_{x_1}(f_{n+1})$\, are contained in \,$\textrm{D}$ we
conclude from the inequalities \ref{ineqinv} and Proposition \ref{invariorbit} that
\begin{eqnarray} 
d\big(f^i_{n+1}(x),f^{i+1}_{n+1}(x)\big)\geq r \quad \mbox{and}  
\quad d\big(f^{i}_{n}(y),f^{i+1}_{n}(y)\big)\geq r \nonumber
\end{eqnarray} 
for all \,$i\in\mathbb{N}$\,. It follows from Lemma \ref{curvagama1} that 
\,$d(\gamma_1,\gamma_0)\geq r$\,.
Consequently,   
\,$\textrm{Int(D}_0)-\textrm{D}_1$\, contains a ball of radius \,$r/3$.

Iterating this procedure we obtain an infinite decreasing nested sequence of closed disks 
\,$\big(\textrm{D}_i\big)_{i\in\mathbb{N}}$\, such that
\,$\textrm{D}_{i+1}\subset \textrm{Int(D}_i)$\, and 
\,$\textrm{Int(D}_i)-\textrm{D}_{i+1}$\, contains a ball of radius \,$r/3$\,, contradicting the fact that 
\,$S^2$\, has finite diameter. This proves the desired result.
\end{proof}


\vskip40pt

\section{Proofs of Theorems \ref{teorcentral} and \ref{teore3}}
\vskip2mm


We  are now  ready to prove our results addressed in the first section. 
Let \,$\{\mathcal{V}_k\}_{k\geq1}$\, be the decreasing nested sequence of \,$C^1$ neighborhoods of the identity map
of \,$S^2$\, as defined in Section 5.
It is easy to see
that Theorem \ref{teorcentral} follows immediately from the following result.

\medskip

\begin{thm}
If \,$G\subset \textrm{\rm Diff}^{\hskip1pt 1}(S^2)$\, is a \,$k$-nilpotent subgroup finitely 
generated by elements in \,$\mathcal{V}_k$\, then
$$Fix(G_{(1)}\,,f_1\,,\ldots\,,f_n)\neq\emptyset \ \ , \ \  \forall f_1\,,\ldots\,,f_n\in\mathcal{V}_k\,.$$
\end{thm}

\begin{proof}
We argue by induction on the nilpotency length of the group \,$G$\, and we will follow the same steps as in the
proof of the Main Lemma.
The case \,$k=1$\, reduces to Bonatti's Theorem.
Assume that the theorem is true for \,$l$-nilpotent subgroups with \,$1\leq l\leq k$\,. In order to prove the theorem
for the nilpotency length of \,$G$\, equal to \,$k+1$\,
 we shall use induction on the number \,$n$\, of diffeomorphisms \,$f_1\,,\ldots\,,f_n$\,.

Suppose \,$G$\, is a  \,$(k+1)$-nilpotent subgroup finitely generated by   elements in \,$\mathcal{V}_{k+1}$.
Let \,$f_1\in\mathcal{V}_{k+1}\cap G$. As in the proof of the Main Lemma the subgroup \,$H$\, 
 generated by \,$G_{(1)}$\, and \,$f_{1}$\, is nilpotent  with length of nilpotency at most \,$k$\, and
\,$G_{(1)}$\, has a finite set of generators
 \,$\{g_1,\ldots,g_m\} \subset \mathcal{V}_k$\,.
The induction assumption on \,$k$\, implies that  
$$Fix(H_{(1)},g_1,\ldots,g_m,f_{1})=Fix(G_{(1)},f_1)$$
is nonempty. This proves the case  \,$n=1$.

Assume the theorem is true for some \,$n\geq 1$\,. 
Suppose by contradiction that there are \,$f_{1},\ldots,f_{n},f_{n+1}\in \mathcal{V}_{k+1} \cap G$\, such that
\begin{eqnarray}\label{hipindteorcentral}
Fix(G_{(1)},f_{1},\ldots,f_{n},f_{n+1})=\emptyset.
\end{eqnarray}

From the induction assumption on $n$ it follows that there is a point  
\,${\tilde p}\in S^2$\, such that 
$${\tilde p}\in Fix(G_{(1)},f_{1},\ldots,f_{n})-Fix{(f_{n+1})}\,.$$
Proposition \ref{invariorbit} implies that  there is a point 
$$p\in \overline{\mathcal{O}^+_{\tilde p}(f_{n+1})} \subset Fix(G_{(1)},f_{1},\ldots,f_{n})$$
which is an \,$\omega$-recurrent point for \,$f_{n+1}$\,. Moreover
$$p\in Fix(G_{(1)},f_{1},\ldots,f_{n})-Fix{(f_{n+1})}\,.$$
Now fix the character curve \,$\gamma\subset\gamma^{\hskip1pt p}_{f_{n+1}}$\, 
 and a closed disk  \,$\textrm{D}$\,  enclosed by \,$\gamma$.
Applying the Main Lemma we obtain  
$$Fix(G_{(1)},f_{1},\ldots,f_{n},f_{n+1})\neq\emptyset\,.$$
This contradicts the equality  \ref{hipindteorcentral}. The theorem  is proved.
\end{proof}

\medskip

Let \,$p\in S^2$\, and \,$f\in\mathcal{V}_1$\, be such that \,$f(p)\neq p$\,. Suppose  the orbit 
\,$\mathcal{O}^+_p(f)$\, is finite. In that case, the point \,$p$\, is a nontrivial \,$\omega$-recurrent point for
\,$f$\, and we can consider the character curve \,$\gamma \subset \gamma^{\,p}_{f}$\, given by Proposition
\ref{propbona1}.

In the proof of the next result we are repeating the arguments used in the proofs of  Main Lemma and Theorem 
\ref{teorcentral}.

To prove Theorem \ref{teore3} it suffices to prove the following result.

\medskip

\begin{thm}
Let \,$G\subset\mbox{\rm Diff}^{\hskip1pt 1}(S^2)$\, be a \,$k$-nilpotent subgroup finitely generated by 
elements in \,$\mathcal{V}_k$\,. Suppose there exists a point \,$p\in S^2$\, with finite nontrivial \,$G$-orbit. If 
\,$f\in \mathcal{V}_k \cap G$\, is such that \,$f(p)\neq p$\, then the \,$G$-action has a fixed point in the interior
of each disk enclosed by the character curve \,$\gamma\subset\gamma^{\,p}_f$.
\end{thm}

\begin{proof}
The proof will be by induction on \,$k$.

For \,$k=1$\, the  group \,$G$\, is commutative. To prove this case let \,$p\in S^2$,
\,$f\in \mathcal{V}_1 \cap G$\, and \,$\gamma\subset\gamma^{\,p}_f$ be as in the theorem.
Let \,$f_1,\ldots,f_n\in\mathcal{V}_1$\, be a set of generators of \,$G$\, and
fix a point
\,$q\in Fix(G)$\, given by Bonatti's Theorem.

Let \,$\textrm{D}\subset S^2-\{q\}$\, be the disk enclosed by \,$\gamma$\,. From Proposition
\ref{propbona1}  we know that there exists a point \,${\tilde p}_1\in Fix(f)\cap\textrm{Int(D)}$\,. If \,${\tilde
p}_1$\, is not a fixed point for \,$G$\, then there exists an integer \,$0\leq \lambda <n$\, such that 
$${\tilde p}_1 \in Fix(f_0,\ldots,f_{\lambda}) \cap \textrm{Int(D)} \quad \mbox{and} \quad f_{\lambda+1}({\tilde
p}_1)\neq {\tilde p}_1$$
where \,$f_0$\, denotes the identity map.

Now, we use the \,$f_{\lambda+1}$-invariance of \,$Fix(f_0,\ldots,f_{\lambda})$\, and Corollary
\ref{curvagama2} to obtain
 a point
\,$p_1\in Fix(f_0,\ldots,f_{\lambda})\cap\textrm{Int(D})$\, which is an \,$\omega$-recurrent point for
\,$f_{\lambda+1}$\,.  If \,$f_{\lambda+1}(p_1)\neq p_1$\,,   once again, Proposition \ref{propbona2} implies that
there exists a point 
\,${\tilde p}_2\in Fix(f_0,\ldots,f_{\lambda},f_{\lambda+1})\cap\textrm{Int(D}_1)$\, where \,$\textrm{D}_1 \subset
\textrm{Int(D})$\, is the disk enclosed by the  character curve \,$\gamma_1\subset \gamma^{\,p_1}_{f_{\lambda+1}}$.
Repeating these arguments no more than 
\,$n$\, times we get a point \,${\tilde q}\in Fix(f_1,\ldots,f_n)\cap \textrm{Int(D})$\, proving the theorem for
the case \,$k=1$.

Now, suppose for some \,$k\geq1$\, the theorem is true for all \,$1\leq l\leq k$\,. We will prove it for \,$k+1$\,.

For this let \,$p\in S^2$, \,$f\in \mathcal{V}_{k+1} \cap G$\, and \,$\gamma\subset\gamma^{\,p}_{f}$\, be as in the
theorem. Let \,$f_1,\ldots,f_n \in \mathcal{V}_{k+1}$\, such that \,$G=\langle f_1,\ldots,f_n \rangle$\, is a 
\,$(k+1)$-nilpotent group and 
fix  \,$q\in Fix(G)$\, given by  Theorem \ref{teorcentral}.

Let \,$H=\langle G_{(1)}, f \rangle$. We know that  \,$H$\, is nilpotent with length of nilpotency at most \,$k$. It
is finitely generated by elements in \,$\mathcal{V}_k$\, and the \,$H$-orbit of \,$p$\, is finite and nontrivial.
Thus, from the induction assumption on \,$k$\,  there exists a point \,${\tilde p}_1 \in Fix(G_{(1)}, f)
\cap \textrm{Int(D})$\, where \,$\textrm{D}$\, is given as above.

Repeating the arguments as in the case \,$k=1$\, and using   Proposition \ref{invariorbit} and Corollary
\ref{curvagama2},  and the Main Lemma we get, after no more then \,$n$\, steps, that there exists a point 
\,${\tilde q}\in Fix(f_1,\ldots,f_n) \cap \textrm{Int(D})$\,, and the proof is finished.
\end{proof}


\end{document}